\documentclass[12pt]{article}
\usepackage{a4, amssymb, exscale}
\usepackage{makeidx, graphicx, multicol}

\newtheorem{lemma}{Lemma}
\newtheorem{proposition}{Proposition}
\newtheorem{theorem}{Theorem}

\newcommand{\Real}{\mathbb{R}}
\newcommand{\Comp}{\mathbb{C}}
\newcommand{\id}{\mathrm{id}}
\newcommand{\Cabb}{\mathrm{C}}
\newcommand{\half}{{\textstyle\frac{1}{2}}}

\newcommand{\supp}{\mathrm{supp}}
\newcommand{\ol}[1]{\overline{#1}}
\newcommand{\per}{\mathrm{per}}
\newcommand{\interior}[1]{\mathrm{int}(#1)}
\newcommand{\Var}{\mathrm{Var}}

\newcommand{\PM}{\mathrm{M}}
\newcommand{\Sing}{\mathrm{S}}
\newcommand{\End}{\mathrm{E}}
\newcommand{\Turn}{\mathrm{T}}
\newcommand{\Surn}{\mathrm{S}}
\newcommand{\Red}{\mathrm{V}}

\newcommand{\definition}[1]{\textit{#1}}

\newcommand{\proof}{{\textit{Proof}\enspace}}

\newcommand{\eop}{\ \vbox{\hrule
                       \hbox{\vrule
                             \hskip 6pt
                             \vrule height 6pt width 0pt
                             \vrule}%
                       \hrule}%
                     \vspace{\medskipamount}
                }

\newlength{\graphicthick}
\newlength{\graphicmid}
\newlength{\graphicthin}

\newsavebox{\ttt}
\sbox{\ttt}{}
\begin{document}
\title{A Note on a Theorem of Parry}
\author{Chris Preston}
\date{\small{}}
\maketitle

\thispagestyle{empty}

\setlength{\graphicthick}{0.25mm}
\setlength{\graphicmid}{0.2mm}
\setlength{\graphicthin}{0.1mm}

In \cite{parry66} Parry shows that a topologically transitive continuous piecewise monotone mapping $f$ with 
positive topological entropy $\mathrm{h}(f)$ is conjugate to a uniformly piecewise linear mapping with slope 
$\exp(\mathrm{h}(f))$.
In this note we generalise Parry's result somewhat to what we call the class of essentially transitive mappings.
This generalisation is of some interest in as much as for mappings with one turning point the converse also holds,
i.e., a uniformly piecewise linear mapping $g$ with one turning point and with slope  $\beta > 1$ is essentially 
transitive. (In fact, $g$ is topologically transitive if and only if the slope $\beta$ lies in the interval 
$(\sqrt{2},2]$; if $\beta \in (1,\sqrt{2}]$ then $g$ is only essentially transitive.) The proof of our 
generalisation relies heavily on a result from the kneading theory of Milnor and Thurston \cite{milnorthurston77},
\cite{milnorthurston88}, which states that a continuous piecewise monotone mapping $f$ with positive topological 
entropy $\mathrm{h}(f)$ is semi-conjugate to a uniformly piecewise linear mapping with slope $\exp(\mathrm{h}(f))$.
We show that if $f$ is essentially transitive then this forces the semi-conjugacy to be a conjugacy.

Let  $a ,\, b \in \Real$  with  $a < b$  and put  $I = [a,b]$. The set of continuous mappings which map the 
interval  $I$  back into itself will be denoted by  $\Cabb(I)$. If  $f \in \Cabb(I)$  and  $n \ge 0$  then  $f^n$
will denote the \definition{$n$\,th iterate} of  $f$, i.e.,  $f^n \in \Cabb(I)$  is defined inductively by  
$f^0(x) = x$, $f^1(x) = f(x)$  and  $f^n(x) = f(f^{n-1}(x))$ for each $x \in I$, or, without arguments, by
$f^0 = \id_I$, $f^1 = f$  and  $f^n = f \circ f^{n-1}$ for each $n \ge 2$.

Mappings $g,\, h \in \Cabb(I)$ are said to be \definition{conjugate} if there exists a homeomorphism  
$\psi : I \to I$  (which in the present situation just means a continuous and strictly monotone mapping of  $I$  
onto itself) such that  $\psi \circ g = h\circ \psi$. Let $f \in \Cabb(I)$; a subset $B$ of $I$ is said to be 
\definition{$f$-invariant} if $f(B) \subset B$, and the mapping $f$ is said to be 
\definition{(topologically) transitive} if whenever $F$ is a closed $f$-invariant subset of $I$ then either  
$F = I$  or  the interior $\interior{F}$ of $F$ is empty. There are several other standard conditions which are 
equivalent to that of being transitive; see, for example, Walters \cite{walters82}.

We are here interested in a special class of mappings from  $\Cabb(I)$, namely the piecewise monotone mappings. 
A mapping  $f \in \Cabb(I)$  is \definition{piecewise monotone} if there exists $p \ge 0$  and  
$a = d_0 < d_1 < \cdots < d_p < d_{p+1} = b$  such that  $f$  is strictly monotone on each of the intervals  
$[d_k,d_{k+1}]$, $k = 0,\,\ldots,\, p$. Let  $f$  be piecewise monotone and suppose the minimal choice for
the  $d_k$'s has been made, i.e., so that $f$  is not monotone (or, equivalently, is not injective) on any open 
interval containing $d_k$ for each $1 \le k \le p$; then  $d_1,\,\ldots,\, d_p$  are called the 
\definition{turning points} of  $f$  and the intervals  $[d_k,d_{k+1}]$, $k = 0,\,\ldots,\, p$  the 
\definition{laps} of  $f$.

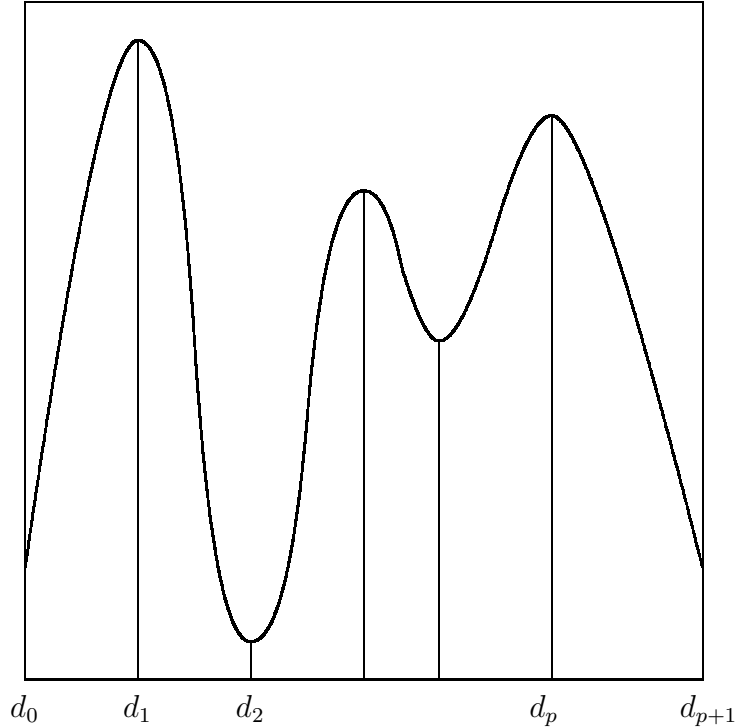
\begin{center}
\setlength{\unitlength}{1.0mm}
\begin{picture}(100,100)
\put(5,5){\framebox(90,90){}}

\linethickness{\graphicthick}

\qbezier(5,20)(15,90)(20,90)
\qbezier(27.5,50)(25,90)(20,90)
\qbezier(27.5,50)(30,10)(35,10)
\qbezier(42.5,40)(40,10)(35,10)
\qbezier(42.5,40)(45,70)(50,70)
\qbezier(55,60)(53,70)(50,70)
\qbezier(55,60)(58,50)(60,50)
\qbezier(67.5,65)(63,50)(60,50)
\qbezier(67.5,65)(72,80)(75,80)
\qbezier(95,20)(80,80)(75,80)

\linethickness{\graphicthin}

\put(20,5){\line(0,1){85}}
\put(35,5){\line(0,1){5}}
\put(50,5){\line(0,1){65}}
\put(60,5){\line(0,1){45}}
\put(75,5){\line(0,1){75}}

\put(3,0){$d_0$}
\put(18,0){$d_1$}
\put(33,0){$d_2$}
\put(72,0){$d_p$}
\put(92,0){$d_{p+1}$}

\end{picture}
\end{center}

The set of piecewise monotone mappings in  $\Cabb(I)$ will be denoted by $\PM(I)$ and for each $f \in \PM(I)$ the 
set of turning points of  $f$ by  $\Turn(f)$. The mappings in $\PM(I)$ are closed under composition: If
$f,\,g \in \PM(I)$ then $g \circ f \in \PM(I)$ and it is easy to see that
$\Turn(g \circ f) = \{ x \in (a,b) : x \in \Turn(f)\ \mbox{or}\ f(x) \in \Turn(g) \}$. In particular, 
$f^n \in \PM(I)$ for all $f \in \PM(I)$, $n \ge 1$. 

Let  $\ell(f)$  denote the number of laps of  $f \in \PM(I)$, so $\ell(f) = \#(T(f)) + 1$; also let  
$\mathrm{h}(f) = \inf_{n\ge 1} n^{-1} \log \ell(f^n)$, and thus $\mathrm{h}(f) \ge 0$.

\begin{lemma}\label{lemma_1}
$\mathrm{h}(f) = \lim\limits_{n\to\infty} n^{-1} \log \ell(f^n)$.
\end{lemma}

\proof
If  $f ,\, g \in \PM(I)$  then  $\ell(f \circ g) \le \ell(f)\ell(g)$, since each lap of $g$  contains at most  
$\ell(f)$  laps of  $f\circ g$, and so in particular $\ell(f^{m+n}) \le \ell(f^m)\ell(f^n)$  for all  
$m,\, n \ge 1$. Put  $a_n = \log \ell(f^n)$; then  $a_{m+n} \le a_m + a_n$  for all  $m,\, n \ge 1$, and hence, as 
is well-known $\lim_{n\to\infty} a_n/n  =  \inf_{n\to\infty} a_n/n$. 
\eop

Misiurewicz and Szlenk  \cite{misiurewiczsz80} show that  $\mathrm{h}(f)$  is the topological entropy of  $f$ 
(which is why we denote this quantity by $\mathrm{h}(f)$).

Let $\beta > 0$;  a mapping $g\in \PM(I)$  is said to be \definition{uniformly piecewise linear} with 
\definition{slope}  $\beta$  if on each of its laps  $g$  is linear with slope  $\beta$  or $-\beta$. The result 
of Parry referred to in the title of this note (Theorem~1 in \cite{parry66}) states that if $f \in \PM(I)$ is 
transitive then $\mathrm{h}(f) > 0$ and $f$ is conjugate to a uniformly piecewise linear mapping with slope 
$\beta = \exp(\mathrm{h}(f))$. Parry's result will be generalised somewhat below to what we call the class of 
essentially transitive mappings in $\PM(I)$.

Let $f \in \Cabb(I)$; a closed set  $C \subset I$  is called an \definition{$f$-cycle}  with \definition{period}
$m \ge 1$  if  $C$  is the disjoint union of non-trivial closed intervals  $B_0,\,\ldots,\, B_{m-1}$  such that  
$f(B_{k-1}) \subset B_k$  for $k = 1,\,\ldots,\, m-1$  and  $f(B_{m-1}) \subset B_0$ (and so in particular $C$ is 
$f$-invariant). An $f$-cycle $C$  is said to be \definition{(topologically) transitive} if whenever $F$  is a 
closed $f$-invariant subset of  $C$ then either  $F = C$  or  $\interior{F} = \varnothing$. Thus $f$ being 
transitive just means that the whole interval $I$ is a transitive $f$-cycle (with period $1$). For each subset 
$B \subset I$ put
\[     \End(B,f) = \{ x \in I : f^n(x) \in B \ \mbox{for some}\  n \ge 0 \}\;,\]
i.e., $\End(B,f)$ consists of those points $x$ for which some iterate $f^n(x)$ lies in the set $B$. The complement
$I \setminus \End(B,f)$ of this set is always $f$-invariant and if $B$ is $f$-invariant then so is $\End(B,f)$.

We call a mapping $f \in \Cabb(I)$ \definition{essentially transitive} if there exists a transitive $f$-cycle $C$ 
such that $I \setminus \End(C,f)$ is countable. 

\begin{theorem}\label{theorem_1}
If the mapping $f \in \PM(I)$ is essentially transitive then $\mathrm{h}(f) > 0$ and $f$ is conjugate to a 
uniformly piecewise linear mapping with slope $\exp(\mathrm{h}(f))$.
\end{theorem}

\proof 
This will follow directly from Theorem~\ref{theorem_2}, Proposition~\ref{prop_2} and Lemma~\ref{lemma_2}. 
Theorem~\ref{theorem_2} is a result from the kneading theory of Milnor and Thurston  \cite{milnorthurston77}, 
\cite{milnorthurston88} and Lemma~\ref{lemma_2} is essentially already a part of Parry's result. Thus the only 
part which is new here is Proposition~\ref{prop_2}. 
\eop

In particular, a transitive mapping is essentially transitive, and so Parry's result is a special case of
Theorem~\ref{theorem_1}. For mappings with one turning point the converse of Theorem~\ref{theorem_1} holds: In 
Proposition~\ref{prop_3} we show that in this case each uniformly piecewise linear mapping with slope $\beta > 1$ 
is essentially transitive.

We should point out that an essentially transitive mapping involves a very special situation, as the follows 
result indicates:

\begin{proposition}\label{prop_1}
Let $f \in \PM(I)$ and let $C$ be an $f$-cycle such that $I \setminus \End(C,f)$ is countable. Then the period of 
$C$ is of the form $2^p$ for some $p \ge 0$. Moreover, there exists $q \ge 0$ such that each periodic point in  
$I \setminus \End(C,f)$  has a period which divides  $2^q$, and  each point in  $I \setminus \End(C,f)$ is 
eventually periodic.
\end{proposition}

\proof
This is Proposition~9.2 in Preston \cite{preston88}, the proof of which is based on an idea occurring in 
Block \cite{block77}, \cite{block79}, Misiurewicz \cite{misiurewicz80}, and in the proof of \v{S}arkovskii's 
theorem (\v{S}arkovskii \cite{sarkovskii64}, \v{S}tefan \cite{stefan77}) given in Block, Guckenheimer, 
Misiurewicz and Young \cite{blockgmy80}.
\eop

Let  $\Red(I) = \{ \psi \in \Cabb(I) : \psi\  \mbox{is increasing and surjective} \}$  (where increasing means 
only that  $\psi(x) \ge \psi(y)$  whenever  $x \ge y$). A pair $(\psi,g)$ with $\psi \in \Red(I)$  and  
$g \in \PM(I)$ is called a \definition{reduction} (or \definition{semi-conjugacy}) of  $f \in \PM(I)$  if  
$\psi\circ f = g\circ\psi$. 

\begin{theorem}\label{theorem_2}
Let $f \in \PM(I)$ with $\mathrm{h}(f) > 0$. Then there exists a reduction  $(\psi,g)$  of  $f$  such that  $g$  
is uniformly piecewise linear with slope  $\exp(\mathrm{h}(f))$.
\end{theorem}

\proof
This can be found in Milnor and Thurston \cite{milnorthurston77}, \cite{milnorthurston88}. A modification of 
their proof (not using any complex analysis) is given at the end of this note. 
\eop

\begin{lemma}\label{lemma_2}
Let $f \in \PM(I)$; if there exists a transitive $f$-cycle then  $\mathrm{h}(f) > 0$. In particular,  
$\mathrm{h}(f) > 0$  whenever  $f$ is  essentially transitive.
\end{lemma}

\proof
Let  $C$  be a  transitive $f$-cycle; put $m = \per(C)$, let  $B$  be one of the  $m$  components of  $C$, and 
let  $g$ be the restriction of  $f^m$  to  $B$, which means that  $g \in \PM(B)$. Then 
$\ell(g^n) \le \ell(f^{mn})$ for each $n \ge 0$ and thus by Lemma~\ref{lemma_1} 
\[      
 \mathrm{h}(f)  =  \lim_{n\to\infty}  \frac{1}{mn} \log \ell(f^{mn}) \ge  
  \lim_{n\to\infty}  \frac{1}{mn} \log \ell(g^{n}) = \frac{1}{m}\, \mathrm{h}(g)\;. 
\]
But $g \in \PM(B)$ is transitive and so it is enough to show that $\mathrm{h}(f) > 0$ for each transitive 
$f \in \PM(I)$ (which is already part of Parry's result). This holds because a transitive mapping $f \in \PM(I)$ 
contains some kind of `horse-shoe': There exist  $p \ge 1$  and  $D,\, E \subset I$ with  
$D \cap E = \varnothing$  such that  $f^p(D) \cap f^p(E) \supset D \cup E$. Thus  $f^p$  is at least $2$  to  $1$  
(i.e., $\#((f^p)^{-1}({x})) \ge 2$ for each  $x \in I$), and so $f^{pn}$  is at least  $2^n$  to  $1$  for each  
$n \ge 1$. It follows that $\ell(f^{pn}) \ge 2^n$  for all  $n \ge 1$, and hence  that 
\[      
\mathrm{h}(f)  =  \lim_{n\to\infty}  \frac{1}{pn} \log \ell(f^{pn}) \ge  \frac{1}{p} \log 2  >  0\;.
\]
The existence of such a `horse-shoe' follows, for example, from Theorem~2.1 in Preston \cite{preston88}, which 
states that a transitive mapping in $\PM(I)$ is either exact or semi-exact. However, the reader can probably
can establish the existence directly without too much trouble. 
\eop

\begin{proposition}\label{prop_2}
If  $(\psi,g)$  is a reduction of an essentially transitive $f \in \PM(I)$ then  $\psi$  must be a homeomorphism, 
and so in particular  $f$  and  $g$  are conjugate.
\end{proposition}

\proof
For each $\psi \in \Red(I)$ put
\begin{eqnarray*} 
\lefteqn{\supp(\psi) = \{ x \in I : \psi(J)\  \mbox{is a non-trivial interval for each open}}\hspace{100pt}\\
                      &&\hspace{120pt}  \mbox{interval  $J \subset I$  containing  $x$} \}\;,
\end{eqnarray*}
which means of course that
\begin{eqnarray*} 
\lefteqn{I \setminus \supp(\psi) 
= \{ x \in I : \ \mbox{there exists an open interval $J \subset I$ containing $x$}}\hspace{100pt}\\
                      &&\hspace{40pt}  \mbox{such that   $\psi(J)$ consists of the single point $\psi(x)$} \}\;,
\end{eqnarray*}
and note that $\psi$ is a homeomorphism if and only if $\supp(\psi) = I$.

\begin{lemma}\label{lemma_3}
For each $\psi \in \Red(I)$ the set $\supp(\psi)$ is non-empty and perfect (i.e., it is closed and contains no 
isolated points).
\end{lemma}

\proof
Clearly $\supp(\psi)$ is closed. Moreover, $\psi$ is constant on each connected component of 
$I \setminus \supp(\psi)$. (Let $J$ be such a component and consider $c,\,d \in J$ with $c < d$. For each 
$x \in [c,d]$ there exists an open interval $J_x$ containing $x$ on which $\psi$ is constant. By compactness 
$[c,d]$ is covered by finitely many of these intervals and hence $\psi(c) = \psi(d)$.) In particular, this means 
that $\supp(\psi)$ is non-empty and perfect (since an isolated point of $\supp(\psi)$ would be the common 
end-point of two connected components of $I \setminus \supp(\psi)$). 
\eop

In fact, each non-empty perfect subset of $I$ is of the form $\supp(\psi)$, i.e., if $D$ is a non-empty perfect
subset of $I$ then there exists a $\psi \in \Red(I)$ with $\supp(\psi) = D$. (This is a classical result in real 
analysis, and can be found, for example, in Carath\'eodory \cite{carath18}. A proof is also given in Preston 
\cite{preston88}, Proposition~11.1.) 

For each $f \in \PM(I)$ put $\Surn(f) = \Turn(f) \cup \{a,b\}$. A subset $A \subset I$ will be called
\definition{$f$-almost-invariant} if $f(A \setminus \Surn(f)) \subset A$.

\begin{lemma}\label{lemma_4}
Let  $(\psi,g)$  be a reduction of a mapping $f \in \PM(I)$. Then $\supp(\psi)$ is $f$-invariant and 
$I \setminus \supp(\psi)$ is $f$-almost-invariant.
\end{lemma}

\proof 
Let $x \in I$ with $f(x) \notin \supp(\psi)$; there thus exists an open interval $J$ containing $f(x)$ such that
$\psi(J)$ consists of the single point $y = \psi(f(x))$. Hence $f^{-1}(J)$ is a neighbourhood of $x$ and so there
exists an open interval $K$ containing $x$ with $K \subset f^{-1}(J)$. Then
$g(\psi(K)) = \psi(f(K)) \subset \psi(J) = \{y\}$, since $f(K) \subset J$ and $g \circ \psi = \psi \circ f$. But 
$\psi(K))$ is connected and $g^{-1}(\{y\})$ is finite, and therefore $\psi(K)$ must consist of the single point 
$\{\psi(x)\}$, i.e.,\ $x \notin \supp(\psi)$. This shows that $\supp(\psi)$ is $f$-invariant.

Now let $x \in (I \setminus \supp(\psi)) \setminus \Surn(f)$; there thus exists an open interval $J$ containing 
$x$ with $J \cap \Surn(f) = \varnothing$ such that $\psi(J)$ consists of the single point $y = \psi(x)$. 
Therefore $f(J)$ is an open interval containing $f(x)$ (since $J \cap \Surn(f) = \varnothing$). But 
$\psi(f(J)) = g(\psi(J)) = g(\{y\})$, since $g \circ \psi = \psi \circ f$, and so $\psi(f(J))$ consists of the 
single point $\psi(\{x\})$, i.e., $f(x) \in I \setminus \supp(\psi)$. This shows that $I \setminus \supp(\psi)$ is 
$f$-almost-invariant. 
\eop

For each $f \in \PM(I)$ let $\mathcal{D}(f)$  denote the set of those non-empty perfect subsets of  $I$  which are 
both $f$-invariant and have an $f$-almost-invariant complement. If $(\psi,g)$ is a reduction of $f$ then by 
Lemmas \ref{lemma_3} and \ref{lemma_4} $\supp(\psi) \in \mathcal{D}(f)$. In fact, the converse also holds: For 
each $\psi \in \Red(I)$ with $\supp(\psi) \in \mathcal{D}(f)$ there exists a unique $g \in \PM(I)$ such that 
$(\psi,g)$ is a reduction of $f$. (This is part of Theorem~5.1 in Preston \cite{preston88}.)

If $C$ is an $f$-cycle then let $C^o$ denote the set obtained by removing the two end-points from each component 
of $C$, so if $m$ is the period of $C$ then $\partial C = C \setminus C^o$ consists of exactly $2m$ points.
The set $C^o$ is not necessarily $f$-invariant but it is easy to see that it is $f$-almost-invariant, which in 
turn easily implies that the open set $\End(C^o,f)$ is $f$-almost-invariant. Of course, 
$\End(C^o,f) \subset \End(C,f)$, since $C^o \subset C$; moreover, $\End(C,f) \setminus E(C^o,f)$ is countable, 
since it is a subset of the countable set $\{ x \in I : f^n(x) \in \partial C \ \mbox{for some}\ n \ge 0 \}$.
Thus $I \setminus \End(C,f)$ is countable if and only if $I \setminus \End(C^o,f)$ is.

\begin{lemma}\label{lemma_5}
Let $f \in \PM(I)$,  $D \in \mathcal{D}(f)$  and  $C$  be a transitive $f$-cycle. Then either  
$\End(C^o,f) \subset D$  or  $\End(C^o,f) \cap D = \varnothing$.
\end{lemma}

\proof
Let  $U \subset I$  be open and $f$-almost-invariant. If  $J$  is a (maximal connected) component of $U$ then 
$f(J \setminus \Surn(f)) \subset U$ and it is easily checked that $f(J \setminus \Surn(f))$ is connected (and in 
fact an open interval). There thus exists a unique component $K$ of $U$ such that $f(\ol{J}) \subset \ol{K}$. 
Iterating this then gives us that for each $n \ge 1$ there exists a unique component $K$ of $U$ such that 
$f^n(\ol{J}) \subset \ol{K}$. A component  $J$  of $U$  is called \definition{periodic} if  
$f^m(\ol{J}) \subset \ol{J}$  for some  $m \ge 1$; the smallest such $m \ge 1$  is called the \definition{period}
of  $J$. A component  $K$  of  $U$  is called \definition{eventually periodic} if  $f^n(\ol{K}) \subset \ol{J}$  
for some periodic component  $J$  of $U$  and some  $n \ge 0$.

Now suppose  $\End(C^o,f) \not\subset D$; then  $\End(C^o,f) \cap (I \setminus D) \ne \varnothing$ and it easily 
follows that $U = C^o \cap (I \setminus D) \ne \varnothing$, so $U$  is a non-empty $f$-almost-invariant open 
subset of $I$. Let $J$ be a component of $U$ which is not eventually periodic and for each $n \ge 0$ let $J_n$ be 
the component of $U$ with $f^n(\ol{J}) \subset \ol{J_n}$. Then the intervals $\{J_n\}_{n\ge 0}$ are disjoint and 
$F = \ol{\bigcup_{n\ge 1} J_n}$ is an $f$-invariant closed subset of $C$ (since $f(\ol{J_n}) \subset \ol{J_{n+1}}$ 
for each $n \ge 1$). But $\interior{F} \ne \varnothing$ (since $J_1 \subset F$) and $F \ne C$ (since 
$J \subset C \setminus F$) and this contradicts the fact that $C$ is transitive. Therefore each component of $U$ 
is eventually periodic and in particular $U$ contains a periodic component. Thus let $J$ be a periodic component 
of $U$ with period $m$; then  $K =  \bigcup_{k=0}^{m-1}  \ol{f^j(J)}$  is a closed $f$-invariant subset of  $C$  
with  $\interior{K} \ne \varnothing$, and hence  $K = C$. But  $C^o \cap D \subset K \setminus U$  and  
$K \setminus U$  is finite. Therefore $C^o \cap D = \varnothing$, because  $D$  is perfect. This implies that  
$\End(C^o,f) \cap D = \varnothing$, since  $D$  is $f$-invariant.  
\eop

We can now complete the proof of Proposition~\ref{prop_2}. Let $C$ be transitive $f$-cycle such that  
$I \setminus \End(C^o,f)$ is countable. By Lemmas \ref{lemma_3} and \ref{lemma_4} $\supp(\psi) \in \mathcal{D}(f)$
and so by Lemma~\ref{lemma_5} either $\End(C^o,f) \cap \supp(\psi)$ is empty or $\End(C^o,f) \subset \supp(\psi)$.
But if  $\End(C^o,f) \cap \supp(\psi) = \varnothing$ then $\supp(\psi)$ is countable, which is not possible
since by the Baire category theorem a non-empty countable closed subset of $I$ must contain an isolated point.
Hence $\End(C^o,f) \subset \supp(\psi)$ and therefore $\supp(\psi) = I$, since $\End(C^o,f)$ is dense in $I$. This
implies that $\psi$ is a homeomorphism. 
\eop

\textit{Proof of Theorem~\ref{theorem_1}:\ }
Let $f \in \PM(I)$ be essentially transitive. Then by Lemma~\ref{lemma_2} $\mathrm{h}(f) > 0$ and hence by 
Theorem~\ref{theorem_1} there exists a reduction  $(\psi,g)$  of  $f$  such that  $g$  is uniformly piecewise 
linear with slope  $\exp(\mathrm{h}(f))$. But by Proposition~\ref{prop_2} $\psi$ is a homeomorphism and therefore 
$f$ and $g$ are conjugate.
\eop

We next note a result of Misiurewicz and Szlenk \cite{misiurewiczsz80} which provides us with an alternative method
of calculating  $\mathrm{h}(f)$  for a mapping  $f \in \PM(I)$. This will show in particular that if  
$g \in \PM(I)$  is uniformly piecewise linear with slope  $\beta > 1$, then  $\mathrm{h}(g) = \log \beta$.
For  $f \in \Cabb(I)$  let
\[ 
\Var(f) = \sup \Bigl\{ \sum_{k=0}^{n-1} |f(x_{k+1})-f(x_k)| :  a = a_0 < x_1 < \cdots < x_n = b \Bigr\}\;. 
\]
If  $f \in \PM(I)$ and $a = d_0 < d_1 < \cdots < d_N < d_{N+1} = b$, where $d_1,\,\ldots,\, d_N$ are the turning 
points of $f$ then clearly
\[
\Var(f) =  \sum_{k=0}^N  |f(d_{k+1})-f(d_k)|\;.
\]
In particular, if  $g \in \PM(I)$  is uniformly piecewise linear with slope $\beta > 0$ then  
$\Var(g) = (b-a)\beta$.

\begin{theorem}\label{theorem_3}
Let  $f \in \PM(I)$; then $\mathrm{h}(f) > 0$ holds if and only if  
\[
\limsup_{n\to\infty} n^{-1} \log \Var(f^n) > 0\;.
\]
Moreover, if $\mathrm{h}(f) > 0$  then
\[ 
\mathrm{h}(f) = \lim_{n\to\infty} n^{-1} \log \Var(f^n)\;.
\]
\end{theorem}

\proof  
This is given in Misiurewicz and Szlenk \cite{misiurewiczsz80}. 
\eop

Let  $g$  be uniformly piecewise linear with slope $\beta > 0$. Then  $\Var(g^n) = (b-a)\beta^n$  for each  
$n \ge 1$, since  $g^n$  is uniformly piecewise linear with slope  $\beta^n$, and thus
$\lim_{n\to\infty} n^{-1} \log \Var(g^n) = \log \beta$. Hence by Theorem~\ref{theorem_3} 
$\mathrm{h}(g) = \log \beta$, provided  $\beta > 1$.

We now consider the case of a mapping $f \in \PM(I)$ having a single turning point, and without loss of 
generality it can be assumed that $f$  takes on its maximum there. Moreover, it will be convenient to assume also 
that  $f(\{a,b\}) \subset \{a,b\}$, i.e., that $f(a) = f(b) = a$. Note that it is always possible to reduce things 
to this case by extending the domain of definition of  $f$  to a larger interval. Moreover, this can be done in 
such a way that for each $x \in (a',a) \cup (b,b')$ the iterates of $x$ end up in $[a,b]$ after finitely many 
steps.

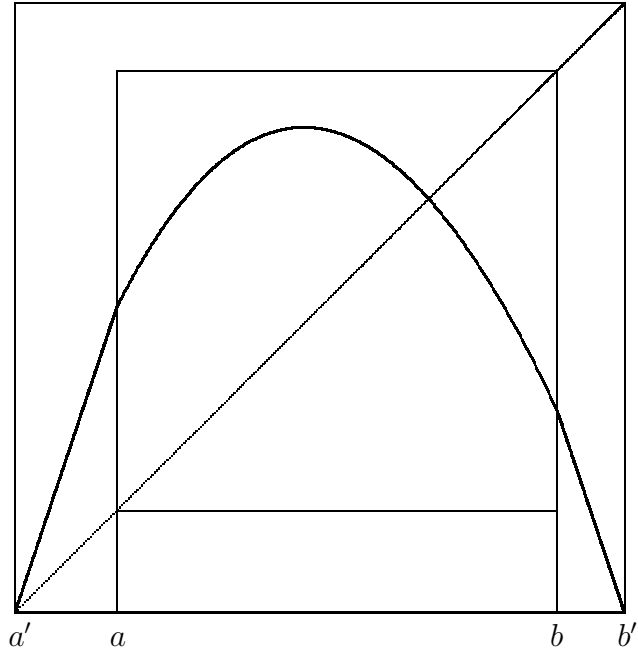
\begin{center}
\setlength{\unitlength}{0.9mm}
\begin{picture}(100,100)
\put(5,5){\framebox(90,90){}}
\qbezier(5,5)(95,95)(95,95)

\linethickness{\graphicthick}

\qbezier(5,5)(20,50)(20,50)
\qbezier(85,35)(95,5)(95,5)
\qbezier(20,50)(50,110)(85,35)

\linethickness{\graphicthin}

\put(20,5){\line(0,1){80}}
\put(85,5){\line(0,1){80}}
\put(20,20){\line(1,0){65}}
\put(20,85){\line(1,0){65}}

\put(4,0){$a'$}
\put(19,0){$a$}
\put(84,0){$b$}
\put(94,0){$b'$}

\end{picture}
\end{center}

Finally, again without loss of generality, assume that $I = [0,1]$; let  $\Sing$  denote the set of mappings  
$f \in \PM([0,1])$ having exactly one turning point and for which  $f(0) = f(1) = 0$.

For each $\beta \in (0,2]$ there is exactly one mapping in $\Sing$ which is uniformly piecewise linear with slope 
$\beta$. This is the mapping $u_\beta$ defined by
\[ 
  u_\beta(x) = \left\{ \begin{array}{cl}
                         \beta x     &\    \mbox{if}\  0 \le x \le \half\;,\\
                        \beta - \beta x  &\      \mbox{if}\  \half \le x \le 1\;. 
\end{array} \right. 
\]

\begin{center}
\setlength{\unitlength}{0.9mm}
\begin{picture}(100,100)
\put(5,5){\framebox(90,90){}}

\linethickness{\graphicthick}

\qbezier(5,5)(50,80)(50,80)
\qbezier(50,80)(95,5)(95,5)

\linethickness{\graphicthin}

\put(50,5){\line(0,1){90}}

\put(5,0){$0$}
\put(49,0){$\half$}
\put(94,0){$1$}

\end{picture}
\end{center}

By Theorem~\ref{theorem_3} $\mathrm{h}(u_\beta) = \log \beta$  for each  $\beta \in (1,2]$. In particular, if  
$\alpha,\, \beta \in (1,2]$  with  $\alpha \ne \beta$  then  $u_\alpha$  and  $u_\beta$  are not conjugate. (If  
$f,\, g \in \PM(I)$  are conjugate then  $\ell(f^n) = \ell(g^n)$  for all  $n \ge 1$, and so  
$\mathrm{h}(f) = \mathrm{h}(g)$.)

Note if $(\psi,g)$ is a reduction of $f \in \Sing$ (with $g \in \PM([0,1])$) then in fact $g \in \Sing$. Thus if 
$f \in \Sing$ with $\mathrm{h}(f) > 0$ then by Theorem~\ref{theorem_2} there exists $\psi \in \Red([0,1])$ such 
that $\psi \circ f = u_\beta \circ \psi$, where $\beta = \log \mathrm{h}(f)$. Moreover, if $f$ is essentially 
transitive then by Theorem~\ref{theorem_1} $f$ and $u_\beta$ are conjugate. In fact, the following result shows 
that the converse is true (since if $f$ and $g$ are conjugate and $f$ is essentially transitive then so is $g$).

\begin{proposition}\label{prop_3}
The mapping $u_\beta$  is essentially transitive for each  $\beta \in (1,2]$.
\end{proposition}

\proof  
Suppose first that  $\beta \in (\sqrt{2},2]$. Let  $C$  be a $u_\beta$-cycle with period  $m$, let  $B$  be one of 
the  $m$  components of  $C$  and let  $g$ be the restriction of  $f^m$  to  $B$. Then  $g$  is uniformly piecewise
linear with slope  $\beta^m$ and, since  $g$  has only one turning point, $\beta^m \le 2$. This is only possible 
if  $m = 1$, i.e. there are no $u_\beta$-cycles with period  $m > 1$. Now let  
$J = [u_\beta^2(\half),u_\beta(\half)]$; then $u_\beta(J) = J$, and so  $J$  is a $u_\beta$-cycle with period  $1$. 
Moreover, if  $K$ is any $u_\beta$-cycle with period  $1$  then  $J \subset K$, since  $\half \in K$. Thus if
$K$  is any $u_\beta$-cycle with  $K \subset J$  then  $K = J$, and from this it is straightforward to show that
$J$  is transitive. (If $J$ is not transitive then there exists a closed $u_\beta$-invariant subset $F$ of $J$ 
with $\interior{F} \ne \varnothing$ and $F \ne J$. Then $U = \interior{F}$ is a non-empty 
$u_\beta$-almost-invariant subset of $J$. An argument similar to that employed in the proof of Lemma~\ref{lemma_5}
shows that $J$ contains a periodic component which can be used to define a $u_\beta$-cycle $K \subset F$, and this 
contradicts the fact that if $K$  is any $u_\beta$-cycle with  $K \subset J$  then  $K = J$.) But it is clear that
$[0,1] \setminus \End(J^o,u_\beta) \subset \{0,1\}$, which shows that $u_\beta$  is essentially transitive. Now 
suppose that  $\beta \in (1,\sqrt{2}]$; then  $d = (1+\beta)^{-1}\beta$  is the unique fixed point of  $u_\beta$
in  $(\half,1)$. Let  $c = 1-d$; thus $u_\beta(c) = d$. Then  $u_\beta^2([c,d]) \subset [c,d]$, and it is easy to 
see that the restriction of  $u_\beta^2$  to  $[c,d]$  is conjugate to  $u_{\beta^2}$.

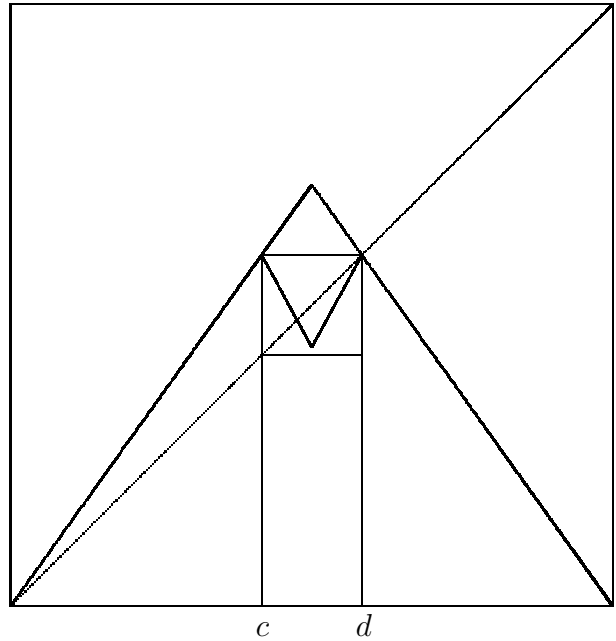
\begin{center}
\setlength{\unitlength}{0.8mm}
\begin{picture}(100,105)
\put(0,0){\framebox(100,100){}}
\qbezier(0,0)(100,100)(100,100)

\linethickness{\graphicthick}

\qbezier(0,0)(50,70)(50,70)
\qbezier(50,70)(100,0)(100,0)

\linethickness{\graphicmid}
\qbezier(41.7,58.3)(50,43)(50,43)
\qbezier(58.3,58.3)(50,43)(50,43)

\linethickness{\graphicthin}

\put(41.7,0){\line(0,1){58.3}}
\put(58.3,0){\line(0,1){58.3}}
\put(41.7,41.7){\line(1,0){16.7}}
\put(41.7,58.3){\line(1,0){16.7}}

\put(40.7,-5){$c$}
\put(57.3,-5){$d$}

\end{picture}

\bigskip
\medskip
\end{center}

But  $u_\beta^m(x) \in [c,d]$  for some  $m \ge 0$  for each  $x \in (0,1)$, and for each  $x \in [0,1]$  the set 
$\{ y \in [0,1] : u_\beta^n(y) = x\  \mbox{for some}\  n \ne 0 \}$  is countable; it thus follows that if  
$u_{\beta^2}$  is essentially transitive then so is $u_\beta$. Therefore  $u_\beta$  is essentially transitive for 
each  $\beta \in (1,2]$.  
\eop

If $\beta \in (0,1]$ then $u_\beta$ is certainly not essentially transitive: In this case it is easy to see that
$Z_*(u_\beta) = [0,1]$, and so there is no transitive $u_\beta$-cycle.

Let $\beta \in (1,2]$; then $u_\beta$ is essentially transitive and thus there exists a transitive $u_\beta$-cycle
$C$ such that $[0,1] \setminus \End(C^o,u_\beta)$ is countable. The proof of Proposition~\ref{prop_3} in fact 
shows that $C$ has period $2^{p}$, where $p \ge 0$ is the smallest integer such that 
$2^{p + 1} \log \beta > \log 2$. Thus $C$ has period $1$ if $\beta \in (\sqrt{2},2]$, period $2$ if 
$\beta \in (\sqrt[4]{2},\sqrt{2}]$, period $4$ if $\beta \in (\sqrt[8]{2},\sqrt[4]{2}]$ and so on. The same then 
holds true for an essentially transitive mapping $f \in \Sing$.

We end this note by giving a proof of Theorem~\ref{theorem_2}. The proof is essentially that to be found in
Milnor and Thurston \cite{milnorthurston77}, \cite{milnorthurston88} but without using any complex analysis.
Fix  a mapping $f \in \PM(I)$  with  $\mathrm{h}(f) > 0$  and put $r = \exp(-\mathrm{h}(f))$; thus  
$r = 1/\beta$  and  $0 < r < 1$. By Lemma~\ref{lemma_1}   $\beta = \lim_{n\to\infty} \ell(f^n)^{1/n}$, and hence  
$r$  is the radius of convergence of the power series  $ \sum_{n\ge 0}  \ell(f^n)t^n$; in particular this means 
that the series  $L(t) =  \sum_{n\ge 0}  \ell(f^n)t^n$  converges for all  $t \in (0,r)$.

\begin{lemma}\label{lemma_101}
$\lim\limits_{t\uparrow r} L(t) = \infty$.
\end{lemma}

\proof
By definition $\ell(f^n) \ge (\exp(\mathrm{h}(f)))^n = \beta^n$  for each  $n \ge 0$ and it therefore follows that
$L(t) \ge  \sum_{n\ge 0}  (\beta t)^n = r(r-t)^{-1}$  for all  $t \in (0,r)$.  
\eop

Let  $\mathcal{J}$  denote the set of non-trivial closed intervals  $J \subset I$; for  $J \in \mathcal{J}$  and  
$n \ge 0$  denote by  $\ell(f^n|J)$  the number of laps of $f^n$  which intersect the interior of  $J$  (and so 
in fact $\ell(f^n|J) = \#(\Turn(f^n) \cap \interior{J}) + 1$). Then  $\ell(f^n|J) \le \ell(f^n|I) = \ell(f^n)$,
and thus in particular the series $L(J,t) =  \sum_{n\ge 0}  \ell(f^n|J)t^n$  converges for all  $t \in (0,r)$. Now,
noting that $L(I,t) = L(t) \ne 0$, put $\Lambda(J,t) = L(J,t)/L(I,t)$  for each  $t \in (0,r)$, thus
$0 \le \Lambda(J,t) \le 1$, because  $L(J,t) \le L(I,t)$.

\begin{lemma}\label{lemma_102}
Let  $J,\, K \in \mathcal{J}$  intersect in a single point. Then
\[             
 \lim_{t \uparrow r} | \Lambda(J \cup K,t) - \Lambda(J,t) - \Lambda(K,t) | = 0\;.
\]
\end{lemma}

\proof
For each  $n \ge 0$
\[       
\ell(f^n|J) + \ell(f^n|K) - 1 \le \ell(f^n|J \cup K) \le \ell(f^n|J) + \ell(f^n|K)\;,
\]
and thus
\begin{eqnarray*}     
| \Lambda(J \cup K,t) - \Lambda(J,t) - \Lambda(K,t) |
         &=& L(I,t)^{-1} | L(J \cup K,t) - L(J,t) - L(K,t) |\\
         &\le& L(I,t)^{-1}  \sum_{n\ge 0}  t^n  = (L(t)(1-t))^{-1}\;.
\end{eqnarray*}
But by Lemma~\ref{lemma_101}   $\lim\limits_{t\uparrow r} (L(t)(1-t))^{-1} = 0$  (since  $r < 1$). 
\eop

\begin{lemma}\label{lemma_103}
Let  $J \in \mathcal{J}$  be such that  $f$  is monotone on  $J$. Then
\[               
\lim_{t\uparrow r} | r\Lambda(f(J),t) - \Lambda(J,t) | = 0\;.
\]
\end{lemma}

\proof
Since  $f$  is monotone on  $f$  it follows that  $\ell(f^{n+1}|J) = \ell(f^n|f(J))$ for each  $n \ge 0$, and thus
\begin{eqnarray*}
 L(J,f) =  \sum_{n\ge 0}  \ell(f^n|J)t^n &=& 1 +  \sum_{n\ge 0}  \ell(f^{n+1}|J)t^{n+1}\\
                   &=& 1 + t  \sum_{n\ge 0}  \ell(f^n|f(J))t^n = 1 + tL(f(J),t)\;.
\end{eqnarray*}
Hence
\begin{eqnarray*}  
\lefteqn{| r\Lambda(f(J),t) - \Lambda(J,t)| 
     = L(I,t)^{-1} | rL(f(J),t) - L(J,t) |}\hspace{50pt}\\
        &\le&  L(I,t)^{-1} (|rL(f(J),t)-tL(f(J),t)| + |tL(f(J),t)-L(J,t)|)\\
        &\le& |r-t| + L(I,t)^{-1} = |r-t| + L(t)^{-1}\;,
\end{eqnarray*}
and by Lemma~\ref{lemma_101}  $\lim_{t\uparrow r} (|r-t| + L(t)^{-1}) = 0$.  
\eop

\begin{lemma}\label{lemma_104}
Let  $J \in \mathcal{J}$. If  $f^m$  is monotone on  $J$  then $\limsup\limits_{t \uparrow r} \Lambda(J,t) \le r^m$.
\end{lemma}

\proof
Since $f$ is monotone on  $f^k(J)$ for  each $k = 0,\,\ldots,\, m-1$  it follows from Lemma~\ref{lemma_103}  that
\[ 
\lim_{t \uparrow r} | r^{k+1}\Lambda(f^{k+1}(J),t) - r^k\Lambda(f^k(J),t) |
             = \lim_{t\uparrow r} r^k | r\Lambda(f(f^k(J)),t) - \Lambda(f^k(J),t) | = 0\;.
\]
Hence  $\lim\limits_{t\uparrow r} | r^m\Lambda(f^m(J),t) - \Lambda(J,t) | = 0$, and so
\[           
\limsup_{t\uparrow r} \Lambda(J,t) = \limsup_{t\uparrow r} r^m \Lambda(f^m(J),t) \le r^m\;,
\]
since  $\Lambda(f^m(J),t) \le 1$  for all  $t \in (0,r)$.  
\eop

\begin{lemma}\label{lemma_105}
There exists a sequence  $\{t_n\}_{n\ge 1}$  from  $(0,r)$  with $\lim_{n\to\infty} t_n = r$  such that  
$\{\Lambda(J,t_n)\}_{n\ge 1}$  converges for all  $J \in \mathcal{J}$.
\end{lemma}

\proof  
Let  $I_o$  be a countable dense subset of  $I$  with  $\{a,b\} \subset I_o$ and  $\Turn(f^n) \subset I_o$  for 
each  $n \ge 1$; let  $\mathcal{J}_o$  be the set of intervals $J = [c,d] \in \mathcal{J}$  such that  
$c ,\, d \in I_o$, thus  $\mathcal{J}_o$  is countable. Now if  $J \in \mathcal{J}$  and  $\{s_n\}_{n\ge 1}$  is 
any sequence from  $(0,r)$  then the sequence  $\{\Lambda(J,s_n)\}_{n\ge1}$ is bounded, and so there exists a 
subsequence $\{n_k\}_{k\ge 1}$  such that  $\{\Lambda(J,s_n)\}_{k\ge 1}$  converges. Therefore, since  
$\mathcal{J}_o$ is countable, a sequence  $\{t_n\}_{n\ge 1}$  from  $(0,r)$  with  $\lim_{n\to\infty}  t_n = r$  
can be found (using the standard diagonal argument) such that $\{\Lambda(J,t_n)\}_{n\ge 1}$  converges for every  
$J \in \mathcal{J}_o$. In fact this sequence then converges for all  $J \in \mathcal{J}$: First consider  
$J = [c,d] \in \mathcal{J}$ with  $c \in I_o$  and  $d \notin I_o$. Let  $\varepsilon > 0$  and choose  $m \ge 1$
so that $r^m < \varepsilon$; since  $d \notin \Turn(f^m) \cup \{a,b\}$  there exist  $u ,\, v \in I_o$  such that
$c < u < d < v$  and  $f^m$  is monotone on  $[u,v]$. Then  $[c,u]$  and $[c,v]$  are both in  $\mathcal{J}_o$, 
and  $\Lambda([c,u],t) \le \Lambda(J,t) \le \Lambda([c,v],t)$  for all  $t \in (0,r)$; hence by Lemmas 
\ref{lemma_102} and \ref{lemma_104} 
\begin{eqnarray*}
 \limsup_{n\to\infty} \Lambda(J,t_n) - \liminf_{n\to\infty} \Lambda(J,t_n)
           &\le& \lim_{n\to\infty} \Lambda([c,v],t_n) - \lim_{n\to\infty} \Lambda([c,u],t_n)\\
           &=& \lim_{n\to\infty} \Lambda([u,v],t_n) \le r^m < \varepsilon.
\end{eqnarray*}
Therefore, since  $\varepsilon > 0$  was arbitrary, the sequence  $\{\Lambda(J,t_n)\}_{n\ge 1}$ converges. The 
same argument also gives that this sequence converges when  $J = [c,d]$  with  $c \notin I_o$  and  $d \in I_o$. 
Finally, if  $J = [c,d]$ and  $c \notin I_o$, $d \notin I_o$  then choose  $u \in I_o$  with  $c < u < d$; then
$\{\Lambda([c,u],t_n)\}_{n\ge 1}$  and  $\{\Lambda([u,d],t_n)\}_{n\ge 1}$  both converge, and so by 
Lemma~\ref{lemma_102}   $\{\Lambda(J,t_n)\}_{n\ge 1}$  converges.  
\eop

Now fix a sequence  $\{t_n\}_{n\ge 1}$  as in Lemma~\ref{lemma_105}, and put  
$\Lambda(J) = \lim_{n\to\infty} \Lambda(J,t_n)$ for each $J \in \mathcal{J}$. Then by Lemmas \ref{lemma_102}, 
\ref{lemma_103} and \ref{lemma_104}:

(1)\enskip If  $J,\,K \in \mathcal{J}$  intersect in a single point then
$\Lambda(J \cup K) = \Lambda(J) + \Lambda(K)$.

(2)\enskip  If  $J \in \mathcal{J}$  and  $f$  is monotone on  $J$  then  $r\Lambda(f(J)) = \Lambda(J)$.

(3)\enskip  If  $J \in \mathcal{J}$, $m \ge 1$  and  $f^m$  is monotone on  $J$  then
$\Lambda(J) \le r^m$.

Also define a mapping $\pi : I \to [0,1]$  by letting  $\pi(a) = 0$  and $\pi(x) = \Lambda([a,x])$  for all 
$x \in (a,b]$.

\begin{lemma}\label{lemma_106}
The mapping  $\pi : I \to [0,1]$  is continuous, increasing and surjective.
\end{lemma}

\proof  
It is clear that  $\pi$  is increasing, and if it is continuous then it is surjective, because  $\pi(a) = 0$  and
$\pi(b) = \Lambda(I) = 1$. Let  $x \in I$ and  $\varepsilon > 0$; choose  $m \ge 1$  so that $r^m < \varepsilon$. 
Then there exists $\delta > 0$  such that  $\{ w \in \Turn(f^m) : |\delta -x| < \delta \} \subset \{x\}$. Now if 
$U = I \cap (x-\delta,x+\delta)$  then  $U$  is a neighbourhood of  $x$  in  $I$, and it follows from (1) and (3) 
that  $|\pi(y)-\pi(x)| < \varepsilon$  for all  $y \in U$, since if  $y > x$  (resp.\ $y < x$) then  $f^m$  is 
monotone on  $[x,y]$ (resp.\ on  $[y,x]$). This shows $\pi$  is continuous.  
\eop

\begin{lemma}\label{lemma_107}
There exists a unique mapping  $\alpha : [0,1] \to [0,1]$  with $\pi \circ f = \alpha \circ \pi$.
\end{lemma}

\proof
If  $\alpha : [0,1] \to [0,1]$  is a mapping with  $\pi \circ f = \alpha \circ \pi$  then $\alpha(z) = \pi(f(x))$
whenever  $x \in I$  is such that  $\pi(x) = z$. Conversely, this relation can be used to define a mapping  
$\alpha$  with  $\pi \circ f = \alpha \circ \pi$, provided $\pi(f(x)) = \pi(f(y))$  whenever  $x ,\, y \in I$  are
such that  $\pi(x) = \pi(y)$. Let  $x ,\, y \in I$  with  $x < y$  and  $\pi(x) = \pi(y)$, and consider  $u,\, v$
with  $x \le u < v \le y$  so that  $f$  is monotone on $[u,v]$. Then  $\pi(u) = \pi(v)$  and hence by (1) and (2)
\[         
r\Lambda(f([u,v])) = \Lambda([u,v]) = \Lambda([a,v]) - \Lambda([a,u])
                        = \pi(v) - \pi(u) = 0\;,
\]
i.e.,  $\Lambda(f([u,v])) = 0$; thus  $\pi(f(u)) = \pi(f(v))$, because  $f(u)$  and $f(v)$  are the end-points of 
the interval  $f([u,v])$. But $[u,v]$ can be written as
\[      
[x,y] = [u_0,u_n] = [u_0,u_1] \cup \cdots  \cup [u_{n-1},u_n]
\]
with  $f$  monotone on each of the intervals  $[u_j,u_{j+1}]$, $0 \le j \le n-1$, and therefore
\[      
\pi(f(x)) = \pi(f(u_0)) = \pi(f(u_1)) = \cdots = \pi(f(u_n)) = \pi(f(y))\;.
\]
This shows that  $\alpha$  exists. The uniqueness of  $\alpha$  follows immediately from the fact that  $\pi$  is 
surjective.  
\eop

Let  $\alpha : [0,1] \to [0,1]$  be the unique mapping with  $\pi\circ f = \alpha \circ \pi$.

\begin{lemma}\label{lemma_108}
The mapping  $\alpha$  is uniformly piecewise linear with slope  $\beta$.
\end{lemma}

\proof
Let  $[c,d]$  be a lap of  $f$  on which  $f$  is increasing, and let  $z \in (\pi(c),\pi(d)]$. Then there exists
$x \in (c,d]$  with  $\pi(x) = z$, and hence by (1) and (2)
\begin{eqnarray*}
      \alpha(z) - \alpha(\pi(c)) &=& \alpha(\pi(x)) - \alpha(\pi(c)) = \pi(f(x)) - \pi(f(c))\\
           &=& \Lambda([a,f(x)]) - \Lambda([a,f(c)]) = \Lambda([f(c),f(x)])
           = \Lambda(f([c,x]))\\ 
          &=& \beta\Lambda([c,x]) = \beta(\Lambda([a,x]) - \Lambda([a,c]))\\
           &=& \beta(\pi(x) - \pi(c)) = \beta(z - \pi(c))\;. 
\end{eqnarray*}
This shows that  $\alpha$  is linear with slope  $\beta$  on  $[\pi(c),\pi(d)]$. If $[u,v]$  is a lap of  $f$  on 
which  $f$  is decreasing then a similar calculation shows that  $\alpha$  is  linear on  $[\pi(u),\pi(v)]$  with 
slope $-\beta$. Therefore  $\alpha$  is uniformly piecewise linear with slope  $\beta$.  
\eop

Now let  $\gamma : [0,1] \to I$  be the linear rescaling given by $\gamma(t) = a + (b-a)t$, and define  
$\psi,\, g : I \to I$  by  $\psi = \gamma \circ\pi$  and $g = \gamma\circ\alpha\circ\gamma^{-1}$. Then  
$\psi \in \Red(I)$  and  $g \in \PM(I)$  is uniformly piecewise linear with slope  $\beta$  (because  $\alpha$  
is); moreover, 
\[              
\psi\circ f = \gamma \circ \pi \circ f 
              = \gamma\circ \alpha\circ \pi = \gamma \circ \alpha\circ \gamma^{-1}\circ \gamma\circ \pi 
              = g\circ \psi\;,
\]
i.e.,  $(\psi,g)$  is a reduction of  $f$. This completes the proof of Theorem~\ref{theorem_2}.  
\eop

\textit{Remark:} In Milnor and Thurston \cite{milnorthurston77}, \cite{milnorthurston88} it is shown that for each
$J \in \mathcal{J}$  there exists a meromorphic function
$L_1(J,\cdot) : D = \{ z \in \Comp : |z| < 1 \}  \to  \Comp \cup \{\infty\}$  with  $L_1(J,t) = L(J,t)$
for all  $t \in (0,r)$. There thus also exists a meromorphic function
$\Lambda_1(J,\cdot) : D  \to  \Comp \cup \{\infty\}$  with  $\Lambda_1(J,t) = \Lambda(J,t)$  for all  $t \in (0,r)$
(with of course $\Lambda_1(J,\cdot) = L_1(J,\cdot)/L_1(I,\cdot)$). Now since  $0 \le \Lambda(J,t) \le 1$  for all
$t \in (0,r)$, it follows that $0 \le \Lambda_1(J,r) \le 1$  and  $\lim_{t\uparrow t} \Lambda(J,t) = \Lambda_1(J,r)$.
In particular, $\Lambda(J) = \Lambda_1(J,r)$. The construction in Lemma~\ref{lemma_105} is therefore not really 
necessary.

\bigskip\bigskip


\bigskip

{\sc Fakult\"at f\"ur Mathematik, Universit\"at Bielefeld}\\
{\sc Postfach 100131, 33501 Bielefeld, Germany}\\
\textit{E-mail address:} \texttt{preston@math.uni-bielefeld.de}\\

\end{document}